\documentclass[10pt, reqno]{amsart}

\usepackage{hyperref} 
\usepackage{amssymb, mathrsfs,bbold}
\usepackage{amsmath, amsthm}
\usepackage{enumerate}
\usepackage{dsfont}
\usepackage{color}
\usepackage[a4paper]{geometry}
\usepackage[utf8]{inputenc}   % To use accentuated letters \UTF{008E}, \UTF{008F}... directly in the Latex file
\usepackage{stmaryrd}
\usepackage{titlesec, wasysym}
\usepackage{appendix}
\usepackage{todonotes, fancyhdr}
\usepackage{chngcntr}
\usepackage{cancel}

\usepackage{tikz-cd}
\usetikzlibrary{calc}
\usepackage{setspace}
\renewcommand{\baselinestretch}{0.99}

\usepackage[all]{xy}
\numberwithin{subsection}{section}
\numberwithin{subsubsection}{subsection}
\numberwithin{equation}{section} % Une equation p.q. est la q-eme equation de la section p.

\allowdisplaybreaks[4]

\newenvironment{Dem}[1][\unskip]{%
    \begin{list}{\hspace{1.15cm}{\sf \textbf{{\small Proof #1 --}}}}{%
        \setlength{\topsep}{0pt}%
        \setlength{\leftmargin}{0pt}%
        \setlength{\rightmargin}{0pt}%
        \setlength{\listparindent}{0pt}%
        \setlength{\itemindent}{0pt}%
        \setlength{\parsep}{0pt}%
        \addtolength{\leftmargin}{20pt}
        \addtolength{\rightmargin}{0pt}%
    } \item }{\hfill $\rhd$\end{list}\smallskip}

\newenvironment{Dem*}[1][\unskip]{%
    \begin{list}{\hspace{0cm}{\sf \textbf{{\small Proof #1 --}}}}{   %
        \setlength{\topsep}{0pt}%
        \setlength{\leftmargin}{0pt}%
        \setlength{\rightmargin}{0pt}%
        \setlength{\listparindent}{0pt}%
        \setlength{\itemindent}{0pt}%
        \setlength{\parsep}{0pt}%
        \addtolength{\leftmargin}{20pt}%
        \addtolength{\rightmargin}{0pt}%
    } \item }{\hfill $\rhd$\end{list}\smallskip}

\renewcommand\thesection       {\arabic{section}}
\renewcommand\thesubsection    {\thesection{\boldmath $.$}\arabic{subsection}}
\renewcommand\thesubsubsection    {\thesection{\boldmath $.$}\arabic{subsection}{\boldmath $.$}\arabic{subsubsection}} %{{\it \thesection}{\boldmath $.$}{\it \arabic{subsection}}{\boldmath $.$}{\it \arabic{subsubsection}}}

\titleformat{\section}[block] %\titleformat{\section}[block]
{\filcenter\normalfont\sffamily\bfseries\Large}  % {\filcenter\normalfont\sffamily\bfseries\Large}
{{\hspace{-0.87cm}}\thesection \hspace{0.2em} --\vspace{0cm}}{0.5em}{} % {{\hspace{-0.87cm}}\thesection \hspace{0.2em} --\vspace{0cm}}{0.5em}{}

\titleformat{\subsection}[runin]
{\filcenter\normalfont\sffamily\bfseries\large}  % {\filright\normalfont\sffamily\bfseries\large}
{{\hspace{0cm}}\thesubsection \hspace{0.15em} -- \vspace{0.1cm}}{.2em}{}   %{\thesubsection \hspace{0.2em} -- \vspace{0.1cm}}{.2em}{}   %{\hspace{-0.7cm}\thesubsection \hspace{0.5em} \vspace{0.3cm}}{.5em}{}  
\titlespacing{\subsection}{-0pc}{1.5ex plus .1ex minus .2ex}{0pc}   %\titlespacing{\subsection}{-0pc}{1.5ex plus .1ex minus .2ex}{0pc}

\titleformat{\subsubsection}[runin]
{\filcenter\normalfont\sffamily\bfseries}   %{\normalfont\sffamily\bfseries}
{\filright\sffamily{\hspace{0cm}}\thesubsubsection\hspace{0.2em} --}{.5em}{}\titlespacing{\subsection}{-0pc}{1.5ex plus .1ex minus .2ex}{0pc}

\newtheoremstyle{mystyle}
{3pt}               %space above
{3pt}               %space below
{\it }                      %bodyfont
{}                      %indent
{\sffamily\bfseries}             %headfont
{}                      %punctuation
{0.5em}                 %space after head
{#1 #2{\hspace{0.2cm}--\hspace{-0.2cm}}  }   %{\llap{#2 }#1{\hspace{0.2cm}--}}

\theoremstyle{mystyle}

\newtheorem{thm}{Theorem}
\newtheorem*{thm*}{Theorem}

\newtheorem{lem}[thm]{\hspace{-0.14cm}  {Lemma} }%[section]
\newtheorem{prop}[thm]{\hspace{-0.13cm} {Proposition}}%[chapter]
\newtheoremstyle{mystyle2}
{3pt}               %space above
{3pt}               %space below
{\it }                      %bodyfont
{}                      %indent
{\sffamily\bfseries}             %headfont
{}                      %punctuation
{0.5em}                 %space after head
{\llap{#2 }#1{\hspace{0.2cm}--}}
\theoremstyle{mystyle2}

\newtheorem*{definition*}{Definition}
\newtheorem*{theorem*}{Theorem}
\newtheorem*{Remark*}{Remark}
\newtheorem*{lem*} {Lemma}
\newtheorem*{defn*} {Definition}
\newtheorem*{prop*} {Proposition}
\newtheorem*{cor*} {Corollary}

\newcommand{\E}{\mathbb{E}}
\newcommand{\N}{\mathbb{N}}
\newcommand{\Z}{\mathbb{Z}}

\newcommand{\R}{\mathbb{R}}
\newcommand{\C}{\mathbb{C}}
\newcommand{\T}{\mathbb{T}}

\renewcommand{\L}{\mathcal{L}}

\newcommand{\F}{\mathcal F}

\newcommand{\CC}{\mathcal C}

\newcommand{\PA}{\mathsf{P}}
\newcommand{\PI}{\mathsf{\Pi}}
\newcommand{\PT}{\overline{\PA}}
\newcommand{\DC}{\mathsf{C}}

\newcommand{\ssk}{\smallskip}

\renewcommand{\epsilon}{\varepsilon}

\newcommand{\bbT}{\textbf{\textsf{T}}}

\newcommand{\mcC}{\mathcal{C}}

\newcommand{\scrS}{\ensuremath{\mathscr{S}}}

\makeatletter
\newcommand*{\defeq}{\mathrel{\rlap{%
                     \raisebox{0.3ex}{$\m@th\cdot$}}%
                     \raisebox{-0.3ex}{$\m@th\cdot$}}%
                     =}
\makeatother

\makeatletter
\newcommand*{\eqdef}{=\mathrel{\rlap{%
                     \raisebox{0.3ex}{$\m@th\cdot$}}%
                     \raisebox{-0.3ex}{$\m@th\cdot$}}%
                     }
\makeatother

\begin{document}

\begin{center}
{\LARGE\sffamily{Non-local quasilinear singular SPDEs  \vspace{0.5cm}}}
\end{center}

\begin{center}
{\sf I. Bailleul} \& {\sf H. Eulry}
\end{center}

\vspace{1cm}

\begin{center}
\begin{minipage}{0.8\textwidth}
\renewcommand\baselinestretch{0.7} \scriptsize \textbf{\textsf{\noindent Abstract.}} We study in this short note a counterpart to the quasilinear generalized parabolic Anderson model (gPAM) on the $2$-dimensional torus where the coefficients are nonlocal functionals of the solution. Under a positivity assumption on the diffusion coefficient we give a local in time solution theory within the framework of paracontrolled calculus.
\end{minipage}
\end{center}

%-----------------------%
\section{Introduction}
\label{SectionIntroNonLocal}
%-----------------------%
We are interested in this work in giving a solution theory for some quasilinear singular stochastic partial differential equations
\begin{equation} \label{EqSPDE}
\partial_t u -A(f(u))\Delta u = B(g(u))\,\xi
\end{equation}
whose coefficients $A(f(u)), B(f(u))$ are {\it nonlocal functionals} of some functions of the solution $u$. We do so here in the mildly singular setting of the $2$-dimensional torus $\T^2$ for a space white noise $\xi$. We need some notations to set the stage of our main result, Theorem \ref{ThmMain} below.

Denote by $(e_1,e_2)$ the canonical basis of $\Z^2\subset\R^2$, and for $i\in\{1,2\}$ and a function $\sigma : \T^2\times\Z^2\rightarrow \C$ set
$$
(D_i\sigma)(x,k) := \sigma(x,k+e_i) - \sigma(x,k),
$$
and ${\bf D}:= (D_1,D_2)$. Write $\mathcal{F}$ for the Fourier transform on $\T^2$ and denote by $\widehat{\sigma}(n,k)=(\mathcal{F}\sigma)(n,k)$ the Fourier transform of the function $\sigma(\cdot,k)$ for each $k\in\Z^2$. Write $\langle k\rangle$ for $1+\vert k\vert$. Given $s\in\R$ and $\alpha>0$ we define the class $\Sigma_\alpha^s$ of functions $\sigma : \T^2\times\Z^2\rightarrow \C$ such that 
\begin{equation} \label{EqFourierBounds}
\big\Vert \partial_x^i {\bf D}^j\sigma(\cdot,k) \big\Vert_{\CC^\alpha} \lesssim \langle k\rangle^{s-\vert j\vert}
\end{equation}
for all multi-indices $i,j$ in $\N^2$, and for which there exists a constant $0<\mu<1$ such that 
\begin{equation} \label{EqSpectralCondition}
(\mathcal{F}\sigma)(n,k) = 0 \quad \textrm{whenever } \quad \vert n\vert >\mu\langle k\rangle.
\end{equation}
For such a symbol, we can easily bound $\widehat{\sigma}$ uniformly in $(n,k)$ by
\begin{align}{\label{EQ:BoundSymbol}}
\big|\widehat{\sigma}(n,k)\big| \lesssim \langle n\rangle^{-4}\,\langle k\rangle^{s}\,\mathbf{1}_{|n|\leq\mu\langle k\rangle}.
\end{align}
The convolution operator with a function $\mathcal{F}^{-1}(\sigma)$ on $\T^2$ whose Fourier transform $\sigma(k)$ satisfies \eqref{EqFourierBounds} provides an elementary example of a function $\sigma\in\Sigma_\alpha^0$ -- here it is independent of $x$ and $\sigma$ is supported on the diagonal $\{n=k\}$. Similarly, ad hoc conditions on a kernel $K$ ensure that the integral operator $v\mapsto \int K(\cdot,y)v(y)dy$ is associated with a symbol $\sigma\in\Sigma_\alpha^s$ for some $s$. A negative $s$ means an operator that regularizes. For $a$ and $b$ in $\Sigma_\alpha^s$ we denote by $A$ and $B$ the pseudodifferential operators with symbols $a$ and $b$ respectively.

We say that $A$ {\it preserves positivity} if $A(v)>0$ whenever $v$ is bounded below by a positive constant. Positivity preserving pseudodifferential operators with a symbol in $\Sigma^0_\alpha$ for some $\alpha>0$ are given by some kernels $\mu(x,\cdot)$ that are finite non-negative measures for all $x\in\T^2$ such that testing $A(v)$ against a smooth test function $\ell$ on $\T^2$ one has
\begin{equation} \label{EqPositivityPreservingOperators}
A(v)(\ell) = \iint v(y)\mu(x,dy)\ell(x)dx.
\end{equation}
This is the case of the integral operator with kernel $K$ when $K$ is bounded below by a positive constant. In the renormalized equation \eqref{EqRenormalizedEquation} below we regularize the space white noise $\xi$ into a smooth function 
$$
\xi^\varepsilon := \mathcal{F}^{-1}\big(\chi(\varepsilon \vert\cdot\vert) \, \mathcal{F}(\xi)(\cdot)\big)
$$
using a `smooth Fourier cut-off' $\chi : [0,\infty)\rightarrow [0,1]$ of class $C^\infty$, equal to $1$ on $[0,1]$ and $0$ on $[2,\infty)$.

\smallskip

\begin{thm} \label{ThmMain}
Pick two regularity exponents $2/3<\beta<\alpha<1$ and $s\leq 0$. Take two symbols $a, b\in\Sigma_\alpha^s$ and assume that the operator $A$ preserves positivity. Take further $f,g\in C^3_b(\R,\R)$.
\begin{enumerate}
    \item[{\it (a)}] There exists two deterministic diverging functions $c_a, c_b : (0,1]\rightarrow C^\infty(\T^2)$ and a positive random time $T$ such that the solution $u^\varepsilon$ to the equation
\begin{equation} \label{EqRenormalizedEquation}
\partial_t u^\varepsilon - A\big(f(u^\varepsilon)\big) \Delta u^\varepsilon = B\big(g(u^\varepsilon)\big) \, \xi^\varepsilon + c_a(\varepsilon) \bigg(\frac{B\big(g(u^\varepsilon)\big)}{A\big(f(u^\varepsilon)\big)}\bigg)^2 f'(u^\varepsilon) - c_b(\varepsilon) \frac{B\big(g(u^\varepsilon)\big)}{A\big(f(u^\varepsilon)\big)} \, g'(u^\varepsilon)
\end{equation}
with initial condition $u_0\in \CC^\alpha(\T^2)$ converges in probability in the parabolic $\alpha$-H\"older space on $[0,T]\times\T^2$, as $\varepsilon\in(0,1]$ goes to $0$.   \vspace{0.1cm}

    \item[{\it (b)}] The singular PDE
    $$
    \partial_t u -A(f(u))\Delta u = B(g(u))\,\xi
    $$
    has a well-defined formulation in a space of paracontrolled functions, where it has a unique solution $(u,u')$. The limit of the $u^\varepsilon$ is given by $u$.
\end{enumerate}
\end{thm}

\smallskip

This statement gives back the result first proved in \cite{OW19, FG19, BDH16} by Otto \& Weber, Furlan \& Gubinelli and Bailleul, Debussche \& Hofmanov\'a when $A$ and $B$ are the identity operators. The approach of Otto \& Weber was later fully developed in \cite{OSSW23, LOT23, LOTT22} to provide a variant of regularity structures well-adapted to a certain class of quasilinear singular stochastic PDEs. It would be interesting to see whether their approach can be twicked to handle nonlocal equations as here. The approach of \cite{BDH16} was developped in a paracontrolled setting in Bailleul \& Mouzard's work \cite{BM23} and fully developped in a regularity structure setting in Bailleul, Hoshino \& Kusuoka's work \cite{BHK23}.

In the proof of Theorem \ref{ThmMain} we obtain point (a) as a corollary of point (b). The proof of point (b) involves two different tasks.
\begin{itemize}
    \item[--] Build an analytical setting where to formulate the equation as a well-defined fixed point problem.
    
    \item[--] Construct a random variable that plays the role of an ill-defined polynomial functional of the noise. It is required as an independent ingredient in the analytical setting.
\end{itemize}

We formulate the equation as a fixed point in a space of paracontrolled functions in Section \ref{SectionPCFormulation}
. We prove therein that it is well-posed locally in time. The construction of the polynomial functional of the noise is done in Section \ref{SectionRenormalization}.

\bigskip

\noindent \textbf{Notations} -- {\it Throughout we denote by $\CC^\gamma(\T^2)$ the Besov space $B_{\infty,\infty}^\gamma(\T^2)$, for any $\gamma\in\R$. For $\alpha\in(0,1)$ and a finite positive time horizon $T$ we set
$$
\mathscr{C}^\alpha_T := C\big([0,T], \CC^\alpha(\T^2)\big) \cap C^{\alpha/2}\big([0,T], L^\infty(\T^2)\big).
$$ 
This space coincides with the $\alpha$-H\"older parabolic space, with an equivalent norm. We will denote by $A$ and $B$ the pseudodifferential operators associated with the symbols $a$ and $b$ respectively. We will use the symbol $\lesssim$ to denote an inequality that holds up to a multiplicative positive constant that does not matter for our purpose. We write $\lesssim_c$ to emphasize that this implicit constant depends only on a parameter $c$.}

%---------------------------------------------------------------%
\section{Paracontrolled formulation of the equation}
\label{SectionPCFormulation}
%---------------------------------------------------------------%

The setting of paracontrolled calculus is now familiar enough that we can go straight to the point. We refer the reader to Gubinelli \& Perkowski's lecture notes \cite{GP17} for an introduction. We refer to the original article \cite{GIP15} of Gubinelli, Imkeller \& Perkowski or the advanced work \cite{BB19} of Bailleul \& Bernicot to see paracontrolled calculus in action in some situations that are more involved than what we need here.

\smallskip

In its basic form a paracontrolled structure is a Banach space of functions of the form
\begin{equation} \label{EqPCStructure}
v = \PA_{v'}X + v^\sharp
\end{equation}
where $v'\in \CC^{\alpha_1}(\T^2)$ with $\alpha_1>0$, and $X\in \CC^\alpha(\T^2)$, and $v^\sharp\in \CC^{\alpha_2}(\T^2)$ with $\alpha_2>\alpha$. Given that $\PA_{v'}X\in \CC^\alpha(\T^2)$ from the continuity estimate of Proposition \ref{PropContinuityParaproduct} in Appendix \ref{APP:PC}, one sees $v^\sharp$ as a `remainder' term in the above decomposition of $v$. This paracontrolled structure is stable by nonlinear functions, a consequence of Bony's paralinearization result, Proposition \ref{Paralinearization} in Appendix \ref{APP:PC}. In the study of the semilinear version of Equation \eqref{EqSPDE}, where $A=1$ say, we can formulate the equation as a fixed point problem in some space of time-dependent paracontrolled functions provided that we have an appropriate definition of the product $X\xi$. An explicit choice of $X$ as a function of $\xi$ brings back the question of that definition to a problem about random variables only. Once this probability problem is understood, solving the equation itself is a purely deterministic problem in a given $\omega$-dependent Banach space of paracontrolled functions, where $\omega$ stands for the chance element. We will follow here a similar strategy, following  Bailleul, Debussche \& Hofmanov\'a's reformulation of the quasilinear equation \eqref{EqSPDE} as a `perturbation' of a semilinear equation \cite{BDH16} 
\begin{equation} \label{EqReformulation}
\partial_tu - A(f(u_0^T))\Delta u = B(g(u))\xi + D_0(u)\Delta u,
\end{equation}
where 
$$
u_0^T := e^{T\Delta}(u_0) \in C^\infty(\T^2)
$$
and 
$$
D_0(u):=A(f(u)) - A(f(u_0^T)).
$$
In order to proceed like that we need first to make sure that both $A(f(v))$ and $B(g(v))$ have a paracontrolled structure if $v$ does and that this structure is preserved by the equation.

\medskip

\subsection{Stability of the paracontrolled structure.}
\hspace{.25cm}Recall the notation of the modified paraproduct $\PT$ from Appendix \ref{APP:PC}. We now fix $X$ for the remainder of this chapter
$$
X= -\Delta^{-1}(\xi),
$$ 
with null mean on $\T^2$. Pick 
$$ 
\frac{2}{3} < \beta<\alpha<1
$$
and define the space $\mathbf{C}^\beta_{\alpha,T}(X)$ as the set of all functions $(u,u')\in\mathscr{C}^\alpha_T \times \mathscr{C}^\beta_T$ such that 
\begin{equation} \label{DEFEQ:ParaStructure}
u^\sharp := \big(u-\PT_{u'}X\big) \in\mathscr{C}^\alpha_T \qquad\text{and}\qquad \sup_{0<t\leq T}t^{\frac{2\beta-\alpha}{2}}\|u^\sharp(t)\|_{\CC^{2\beta}}<+\infty
\end{equation}
which we equip with the norm
$$
\|(u,u')\|_{\mathbf{C}^\beta_{\alpha,T}(X)} := \|u'\|_{\mathscr{C}^\beta_T} + \|u^\sharp\|_{\mathscr{C}^\alpha_T} + \sup_{0<t\leq T}t^{\frac{2\beta-\alpha}{2}}\|u^\sharp(t)\|_{\CC^{2\beta}}.
$$
We will look for a solution to Equation \eqref{EqSPDE}, equivalently Equation \eqref{EqReformulation}, as an element of $\mathbf{C}^\beta_{\alpha,T}(X)$.
\smallskip
The reformulation \eqref{EqReformulation} of Equation \eqref{EqSPDE} involves a number of operations. We need to check how the  a priori paracontrolled structure of a solution to \eqref{EqReformulation} interacts with these operations.

\medskip

\textbf{{a) Pseudodifferential operators and paraproducts.}} The class of symbols $\Sigma^s_\alpha$ is defined from the constraints \eqref{EqFourierBounds} and \eqref{EqSpectralCondition}. For $\alpha_1>0$ and $h_1\in \CC^{\alpha_1}(\T^2)$ we check that the operator $\PA_{h_1}(\cdot)$ has a symbol in the class $\Sigma^0_{\alpha_1}$. We read on the condition \eqref{EqFourierBounds} that the class of operators associated to $\Sigma^s_\alpha$ is a subclass of the set of pseudodifferential operators $\Psi^s_{11}$. The operators in this class send continuously any $\CC^\gamma(\T^2)$ space into $\CC^{\gamma-s}(\T^2)$ if $\gamma>0$ and $(\gamma-s)>0$ are non-integers -- this is the classical Schauder estimate.

One gets as follows the structure \eqref{EqPositivityPreservingOperators} of the operators with a symbol in $\Sigma^0_\alpha$ that are positivity preserving. For $v\in \CC^\gamma(\T^2)$ with $\gamma$ one has from the classical Schauder estimate $A(v)\in \CC^\gamma(\T^2)$. Now for $v\in \CC^\gamma(\T^2)$, since $2 \Vert v\Vert_\infty \pm v \geq \Vert v\Vert_\infty>0$ the positivity preserving property yields $A(2\Vert v\Vert_\infty\pm v)\geq 0$, that is $\pm A(v) \leq 2\Vert v\Vert_\infty A({\bf 1})$, i.e. $\vert A(v)\vert \leq 2\Vert v\Vert_\infty A({\bf 1})$. For any $x\in\T^2$ the linear form $A(\cdot)(x)$ on $\CC^\gamma(\T^2)$ thus has a continuous extension into a linear form on $C^0(\T^2)$. Riesz representation theorem then gives the existence of a non-negative measure $\mu(x,\cdot)$ such that $A(\cdot)(x)$ is the integration map against the measure $\mu(x,\cdot)$.

\smallskip

To check that the paracontrolled structure \eqref{EqPCStructure} is stable by the nonlocal maps $A$ associated with some symbols $\sigma\in\Sigma^s_{\alpha_1}$ we can use the following commutator estimate from Proposition 3.8 in  Bailleul, Dang, Ferdinand \& T\^o's work \cite{BDFT23}, together with the classical Schauder estimate.

\smallskip

\begin{prop}
Pick $s\leq 0, \alpha_1\in (0,1)$ and $\alpha_2 > 0$. For any $\sigma\in\Sigma^s_{\alpha_1}$ and $h_1\in \CC^{\alpha_1}(\T^2)$, for any $h_2\in \CC^{\alpha_2}(\T^2)$ one has
$$
\big\Vert Op(\sigma)(\PA_{h_1} h_2) - \PA_{h_1}\big(Op(\sigma)(h_2)\big) \big\Vert_{\CC^{\alpha_1+\alpha_2-s}} \lesssim_a \Vert h_1\Vert_{\CC^{\alpha_1}} \Vert h_2\Vert_{\CC^{\alpha_2}}.
$$
\end{prop}

\ssk

\textbf{{b) Structural elements in the right hand side of \eqref{EqReformulation}.}} From Proposition \ref{Paralinearization} on paralinearization, Proposition \ref{PropIteratedParaproducts} on iterated paraproducts and the continuity result for the corrector in Proposition \ref{Corrector}, one can expand the first ill-defined product
\begin{align}{\label{EQ:ExpansionB}}
B(f(u)) \, \xi = \PA_{B(f(u))} \xi + \PA_\xi B(f(u)) + u'g'(u) \PI(B(X),\xi)+g'(u)\PI(B(u^\sharp),\xi) + \sharp_B
\end{align}
where the ill-defined product $\PI(B(X),\xi)$ is to be understood as its renormalized counterpart defined in Section \ref{SectionRenormalization} -- note that $u'g'(u) \PI(B(X),\xi)$ is the only term featuring this ill-defined product. The term
\begin{align*}
\sharp_B &= \DC(u'g'(u),B(X),\xi) + \DC(g'(u),u^\sharp,\xi) + \PI\big([B,\PA_{u'g'(u)}](X)\,,\,\xi\big) + \PI\big([B,\PA_{g'(u)}](u^\sharp)\,,\,\xi\big)\\
&\qquad+ \PI\big(B\big(\PA_{g'(u)}(\PT_{u'}X) - \PA_{u'g'(u)}X\big)\,,\,\xi\big) + \PI\big(B\big(g(u)-\PA_{g'(u)}u\big),\xi\big)
\end{align*}
denotes an element of the space $C(\big[0,T],\CC^{2\alpha+\beta-2}(\T^2)\big)$ that is polynomial in $(u,u')$ and depends continuously on $\xi$. Given $\gamma\in\R$ write $=_\gamma$ to mean equality up to an element of $C\big([0,T],\CC^\gamma(\T^2)\big)$. One also has for the Laplace term
\begin{equation*} \begin{split}
    D_0(u)\Delta u&=D_0(u)\big(\Delta\PT_{u'}X+\Delta u^\sharp\big)\\
    &=D_0(u)\big(-\PT_{u'}\xi+[\Delta,\PT_{u'}]X+\Delta u^\sharp\big)\\
    &=_{\alpha+\beta-2} D_0(u)\big(-\PT_{u'}\xi + \Delta u^\sharp\big)\\
    &=_{\alpha+\beta-2} -\PA_{D_0(u)} \big(\PT_{u'}\xi\big) - \PI\big(D_0(u),\PT_{u'}\xi\big) + D_0(u)\Delta u^\sharp   \\
    &=_{\alpha+\beta-2} -\PA_{D_0(u)u'}\xi - \PI\big(A(f(u)),\PT_{u'}\xi\big) + D_0(u)\Delta u^\sharp   \\
    &=_{\alpha+\beta-2} -\PA_{D_0(u)u'}\xi-u'\PI\big(A(f(u)),\xi\big) + D_0(u)\Delta u^\sharp,
\end{split} \end{equation*}
where we used successively the paracontrolled expansion of $u$, the commutation lemma \ref{CommutationParaproducts} to get rid of the $[\Delta, \PT_{u'}]$ term and switch back and forth between $\PT$ and $\PA$, the commutation property \ref{PropIteratedParaproducts} to merge the paraproducts, the fact that, as $D_0(u)=A(f(u))-A(f(u_0^T))$ with $u_0^T$ smooth, so is $A(f(u_0^T))$, and the $\PI\big(A(f(u_0^T)),\PT_{u'}\xi\big)$ term has $C_T\CC^{\alpha+\beta-2}$ regularity, and finally the corrector estimate \ref{Corrector} for the remaining resonant term. The resonant product in the last step can be expanded as follows using the paracontrolled expansion of $u$ and paralinearization lemma \ref{PropRefinedParalinearization}
\begin{align}{\label{EQ:ExpansionA}}
\PI\big(A(f(u)),\xi\big) =_{2\alpha+\beta-2} u'f'(u)\PI\big(A(X),\xi\big) + f'(u)\PI\big(A(u^\sharp),\xi\big).
\end{align}
Note that $-(u')^2f'(u)\PI\big(A(X),\xi\big)$ is the only term in this decomposition of $D_0(u)\Delta u$ featuring the renormalized product $\PI\big(A(X),\xi\big)$. This gives
$$
D_0(u)\Delta u = -\PA_{D_0(u)u'}\xi - u'^2f'(u) \PI\big(A(X),\xi\big) + u'f'(u) \PI\big(A(u^\sharp),\xi\big) + D_0(u)\Delta u^\sharp + \sharp_A
$$
and
\begin{align*}
    \sharp_A &= D_0(u)\,  [\Delta,\PT_{u'}] (X) + \big(\PA_{D_0(u)u'}\xi - D_0(u)\,\PT_{u'}X\big) - \PA_{\PT_{u'}\xi}D_0(u) + \PI\big(\PT_{u'}\xi, A(f(u_0^T))\big)   \\
    &\quad- \big\{\PI\big(\PT_{u'}\xi,A(f(u))\big) - u' \PI\big(\xi,A(f(u))\big)\big\} - u' \, \DC\big(u'f'(u),A(X),\xi\big) + u' \,\DC\big(f'(u),u^\sharp,\xi\big)   \\
    &\quad+ u' \, \PI\big([A,\PA_{u'f'(u)}](X) , \xi\big) + u' \, \PI\big([A,\PA_{f'(u)}](u^\sharp) , \xi\big)   \\
    &\quad+ u' \,\PI\big(A\big(\PA_{f'(u)}\PT_{u'}X - \PA_{u'f'(u)}X\big) , \xi\big) + u' \, \PI\big(A\big(f(u)-\PA_{f'(u)}u\big) , \xi\big)
\end{align*}
is an element of $C\big([0,T],\CC^{2\alpha+\beta-2}(\T^2)\big)$ that is polynomial and continuous in $(u,u')$ and continuous in $\xi$. Note that although $(\Delta u^\sharp)(t)$ is an element of $\CC^{2\beta-2}(\T^2)$ at every positive time $t$ it cannot be considered as a remainder term as its $\CC^{2\beta-2}$-norm blows up when $t>0$ goes to $0$. Using Lemma \ref{CommutationWeightedLaplacian} we can then rewrite Equation \eqref{EqReformulation} as a system of coupled equation

\begin{equation*} \begin{split}
\partial_tu - A(f(u_0^T))\Delta u &= \PA_{B(g(u)) - A(f(u))u'} \xi + \big\{A(f(u)) - A(f(u_0^T))\big\}\Delta u^\sharp   \\
&\quad+ g'(u) \PI\big(B(u^\sharp) , \xi\big) - u'f'(u)\,\PI\big(A(u^\sharp) , \xi\big) + \sharp_1(u,u')
\end{split} \end{equation*}
and
\begin{equation*} \label{EqDefnPhiSharp} \begin{split}
\partial_tu^\sharp - A(f(u_0^T))\Delta u^\sharp &= \PA_{B(g(u)) - A(f(u))u'} \xi + \big\{A(f(u)) - A(f(u_0^T))\big\}\Delta u^\sharp   \\
&\quad+g'(u)\PI\big(B(u^\sharp) , \xi\big) - u'f'(u)\,\PI\big(A(u^\sharp) , \xi\big) + \sharp_2(u,u')
\end{split} \end{equation*}
where $\sharp_1(u,u')$ and $\sharp_2(u,u')$ are some elements of $C\big([0,T], \CC^{\alpha+\beta-2}(\T^2)\big)$ that are non-linear functions of $u$ that are multilinear in $(u',u^\sharp)$ and depend continuously on the enhanced noise data 
$$
\big(\xi,\PI(A(X),\xi),\PI(B(X),\xi)\big).
$$
We emphasized in Equations \eqref{EQ:ExpansionB} and \eqref{EQ:ExpansionA} the dependency on the renormalized products $\PI(B(X),\xi)$ and $\PI(A(X),\xi)$ respectively. This will be useful when we investigate the consistency of the renormalization procedure from a dynamical point of view.

\medskip

\textbf{{c) Fixed point formulation.}} We are now able to formulate equation \eqref{EqReformulation} as a fixed point problem. We define on $\mathbf{C}^\beta_{\alpha,T}(X)$ a map $\Phi$ setting
$$
\Phi\big((u,u')\big) = (v,v')
$$
where
$$
v' := \frac{B(g(u)) - \big\{A(f(u)) - A(f(u_0^T))\big\}u'}{A(f(u_0^T))}
$$
and $v$ solves the equation
\begin{equation} \label{EqEquationU} \begin{split}
\partial_tv - A(f(u_0^T))\Delta v &= \PA_{A(f(u_0^T))v'}\xi + \big\{A(f(u)) - A(f(u_0^T))\big\}\Delta u^\sharp + g'(u) \, \PI\big(B(u^\sharp) , \xi\big)  \\
&\quad- u'f'(u) \, \PI\big(A(u^\sharp) , \xi\big) + \sharp_1(u,u')
\end{split} \end{equation}
with the initial value 
$$
v(0) = u_0.
$$ 
Note that in view of Lemma \ref{LEM:StabilityStructure} below, $(v,v')$ is an element of $\mathbf{C}^\beta_{\alpha, T}(X)$. Furthermore the equation on the remainder $v^\sharp$ writes
\begin{equation} \label{EqEquationRemainder} \begin{split}
\partial_tv^\sharp - A(f(u_0^T))\Delta v^\sharp &= \big\{A(f(u)) - A(f(u_0^T))\big\}\Delta u^\sharp + g'(u) \, \PI\big(B(u^\sharp),\xi\big)  \\
&\quad- u'f'(u) \, \PI\big(A(u^\sharp) , \xi\big) + \sharp_2(u,u')
\end{split} \end{equation}
with initial value 
$$
v^\sharp(0) = u_0 - {\sf P}_{v'(0)}X = u_0 - \PT_{v'(0)}X.
$$ 
Denote by $\mathbf{\mathcal{B}}_T(r)$ the following ball of $\mathbf{C}^\beta_{\alpha,T}(X)$
$$
\mathbf{\mathcal{B}}_T(r) := \left\{(u,u')\in\mathbf{C}^\beta_{\alpha,T}(X)\,; u(0) = u_0\,, u'(0) = \frac{B(g(u_0))}{A(f(u_0))}\,, \|u\|_{\mathbf{C}^\beta_{\alpha,T}(X)}\leq r\right\}.
$$
We prove below that for a sufficiently large radius $r$ and sufficiently small time horizon $T$ the map $\Phi$ is a contraction of $\mathbf{\mathcal{B}}_T(r)$.

\medskip

\subsection{Short time estimate.}
\hspace{.25cm}The technical part of the argument is handled by the following Schauder-like estimate found in Lemma 5 of \cite{BDH16}.

\begin{lem}
\label{SchauderWeight}
Pick $u_0\in \CC^\alpha(\T^2)$ and $w\in C^2(\T^2)$ bounded below by some positive constant. Let also $\phi_1,\phi_2$ be such that 
    $$
    \phi_1\in C\left((0,T],\CC^{2\beta-2}(\T^2)\right)\quad \text{and}\quad \ell_1(T):= \sup_{0<t\leq T}t^{\frac{2\beta-\alpha}{2}}\|\phi_1(t)\|_{\CC^{2\beta-2}} < +\infty,
    $$
    and
    
    $$
    \phi_2\in C\left((0,T],\CC^{\alpha+\beta-2}(\T^2)\right)\quad \text{and}\quad \ell_2(T):= \sup_{0<t\leq T}t^{\frac{2\beta-\alpha}{2}}\|\phi_2(t)\|_{\CC^{\alpha+\beta-2}} < +\infty.
    $$
    Then for $T$ small enough the solution of the equation 
    $$
    \partial_tu-w\Delta u=\phi_1+\phi_2
    $$
    with initial data $u_0$ satisfies the estimate
    \begin{eqnarray*}
        \|u\|_{\mathscr{C}^\alpha}+\sup_{0<t\leq T}t^{\frac{2\beta-\alpha}{2}}\|u(t)\|_{\CC^{2\beta}}\lesssim \|u_0\|_{\CC^\alpha} + \ell_1(T) + T^{\frac{\alpha-\beta}{2}} \ell_2(T),
    \end{eqnarray*}
    where the implicit constant depends only on $\Vert w\Vert_{\CC^\alpha}$.
\end{lem}

\smallskip

We obtain as a consequence of Lemma \ref{SchauderWeight} the stability of the mapping $\Phi$ as a function form the space of paracontrolled functions $\mathbf{C}^\beta_{\alpha,T}(X)$ into itself.

\smallskip

\begin{lem}{\label{LEM:StabilityStructure}}
Let $u_0,w,\phi_1,\phi_2$ be as in Lemma \ref{SchauderWeight} and let $u'\in\CC^\beta$. Denote by $u$ the solution to 
$$
\partial_tu-w\Delta u = \PA_{wu'}\xi+\phi_1 + \phi_2
$$
starting from $u(0)=u_0$. Then $(u,u')\in \mathbf{C}^\beta_{\alpha,T}(X)$.
\end{lem}

\begin{Dem}
Set $u^\#:=u-\PT_{u'}X$, it only remains to check that $u^\sharp$ has the regularity matching the definition \eqref{DEFEQ:ParaStructure} to ensure that $(u,u')$ has indeed a paracontrolled structure. Note that in view of Lemma \ref{CommutationWeightedLaplacian} the term 
$$
(\partial_t - w\Delta)\PT_{u'}X-\PA_{wv'}(-\Delta X)
$$
belongs to $C_T\CC^{\alpha+\beta-2}$. As such $u^\sharp$ satisfies the equation
\begin{align*}
\partial_t u^\sharp-w\Delta u^\sharp &= \partial_t u-w\Delta u - (\partial_t - w\Delta)\big(\PT_{u'}X\big)\\
&=\PA_{wu'}\xi+\phi_1 + \phi_2 - \big((\partial_t - w\Delta)\PT_{u'}X-\PA_{wv'}(\xi)\big) - \PA_{wv'}\xi\\
&=_{\alpha+\beta-2} \phi_1 + \phi_2
\end{align*}
starting from $u^\sharp(0)$. Lemma \ref{SchauderWeight} then finishes the proof.
\end{Dem}

\smallskip

$\bullet$ First to control $v'$ let $(u_1,u_1'),(u_2,u_2')\in \mathbf{\mathcal{B}}_T(r)$, we control the difference $v_1'-v_2'$ writing
$$
A(f(u_0^T))\big(v_1'-v_2'\big) = B(g(u_1)) - B(g(u_2)) - \big\{A(f(u_1)) - A(f(u_0^T))\big\} u_1' + \big\{A(f(u_2)) - A(f(u_0^T))\big\}u_2'.
$$
First use the fact that $B$ gains $-s$ derivatives and Lemma \ref{Paralinearization} to get
\begin{eqnarray*}
    \big\| B(g(u_1)) - B(g(u_2)) \big\|_{\mathscr{C}^{\beta}}&\lesssim& \big\| B(g(u_1)) - B(g(u_2)) \big\|_{\mathscr{C}^{\beta-s}}\\
    &\lesssim&\big(1+\|u_1\|_{\mathscr{C}^{\beta}}\big) \|u_1-u_2\|_{\mathscr{C}^{\beta}}   \\
    &\lesssim&T^{\frac{\alpha-\beta}{2}} \big(1+\|u_1\|_{\mathscr{C}^{\alpha}}\big) \|u_1-u_2\|_{\mathscr{C}^\alpha},
\end{eqnarray*}
where we used the assumption that $\beta-s<\alpha$. Then writing 
$$
A_i := A(f(u_i)) - A(f(u_0^T))
$$ 
we have
\begin{eqnarray*}
    \big\| A_1u_1' - A_2u_2' \big\|_{\mathscr{C}^{\beta}} &\lesssim&\|A_1(u_1'-u_2')\|_{\mathscr{C}^{\beta}}+\|(A_1-A_2)u_2'\|_{\mathscr{C}^{\beta}}   \\
    &\lesssim&\|A_1\|_{C_TL^\infty}\|u_1'-u_2'\|_{\mathscr{C}^\beta}+\|A_1\|_{\mathscr{C}^{\beta}}\|u_1'-u_2'\|_{C_TL^\infty}   \\
    &&\quad+\|A_1-A_2\|_{C_TL^\infty}\|u_2'\|_{\mathscr{C}^\beta}+\|A_1-A_2\|_{\mathscr{C}^{\beta}}\|u_2'\|_{C_TL^\infty}.
\end{eqnarray*}
The difference $A_1-A_2 = A(f(u_1)) - A(f(u_2))$ is dealt with in the same way as $B$ above, with $A_i$ estimated per se as we cannot extract a contracting factor from the difference $\|u'_1-u'_2\|_{\mathscr{C}^\beta}$. We have
$$
A_i = \big\{A(f(u_i)) - A\big(f(u(0))\big)\big\} + \big\{A\big(f(u(0))\big) - A(f(u_0^T))\big\},
$$
thus reasoning as for the $B$ term we can use the continuity properties of $A$ and the paralinearization Lemma \ref{Paralinearization} to obtain
$$
\big\| A(f(u_i)) - A\big(f(u(0))\big) \big\|_{\mathscr{C}^{\beta}}\lesssim T^{\frac{\alpha-\beta}{2}} \big(1+\|u_i\|_{\mathscr{C}^{\alpha}}\big) \|u_i-u(0)\|_{\mathscr{C}^\alpha}.
$$
In a similar way, using the classical estimate \eqref{EqHeatEstimate} on $\text{Id}-e^{T\Delta}$ one gets
$$
\big\| A\big(f(u(0))\big) - A(f(u_0^T)) \big\|_{\CC^{\beta}}\lesssim T^{\frac{\alpha-\beta}{2}} \big(1+\|u_0\|_{\CC^\alpha}\big) \|u_0\|_{\CC^\alpha}.
$$
In the end one has
$$
\big\| A_1u_1'-A_2u_2' \big\|_{\mathscr{C}^{\beta}} \leq T^{\frac{\alpha-\beta}{2}}P\big(\|u_i\|_{\mathscr{C}^\alpha},\|u_i'\|_{\mathscr{C}^\beta},\|u_0\|_{\CC^\alpha}\big) \Big(\|u_1-u_2\|_{\mathscr{C}^\alpha}+\|u'_1-u'_2\|_{\mathscr{C}^\beta}\Big),
$$
for some constant $P(\|u_i\|_{\mathscr{C}^\alpha},\|u_i'\|_{\mathscr{C}^\beta},\|u_0\|_{\alpha})$ that grows at most polynomially as a function of its arguments.

\smallskip

$\bullet$ We now deal with the remainder term $v^\sharp$. Recall from \eqref{EqEquationU} and \eqref{EqEquationRemainder} the notations for the functions $\sharp_1$ and $\sharp_2$. Note that $v_1^\sharp-v_2^\sharp$ is a solution of the equation
\begin{equation*} \begin{split}
\big(\partial_t - A(f(u_0^T))\Delta\big)\big(v_1^\sharp-v_2^\sharp\big) &= \big\{A(f(u_1)) - A(f(u_0^T))\big\} \Delta u_1^\sharp - \big\{A(f(u_2)) - A(f(u_0^T))\big\}\Delta u_2^\sharp   \\
&\quad+ g'(u_1) \, \PI\big(B(u_1^\sharp) , \xi\big) - g'(u_2)\,\PI\big(B(u_2^\sharp) , \xi\big)   \\
&\quad+ u_2'f'(u_2)\,\PI\big(A(u_2^\sharp) , \xi\big) - u_1'f'(u_1) \, \PI\big(A(u_1^\sharp) , \xi\big)   \\
&\quad+ \sharp_2(u_1,u_1') - \sharp_2(u_2,u_2')
\end{split} \end{equation*}
with null initial condition. Set
$$
\phi_1 := \big\{A(f(u_1)) - A(f(u_0^T))\big\} \Delta u_1^\sharp - \big\{A(f(u_2)) - A(f(u_0^T))\big\}\Delta u_2^\sharp
$$
and 

\begin{equation*} \begin{split}
\phi_2 &:= g'(u_1) \, \PI\big(B(u_1^\sharp) , \xi\big) - g'(u_2)\,\PI\big(B(u_2^\sharp) , \xi\big)   \\
&\quad+ u_2'f'(u_2)\,\PI\big(A(u_2^\sharp) , \xi\big) - u_1'f'(u_1) \, \PI\big(A(u_1^\sharp) , \xi\big) + \sharp_2(u_1,u_1') - \sharp_2(u_2,u_2').
\end{split} \end{equation*}
We now check that $\phi_1$ and $\phi_2$ satisfy the regularity assumptions of Lemma \ref{SchauderWeight}. 

-- Concerning $\phi_1(t)$, the estimate in $\CC^{2\beta-2}(\T^2)$ is similar to what we did for $v_1'-v_2'$ so one has
$$
\|\phi_1(t)\|_{\CC^{2\beta-2}}\lesssim T^{\frac{\alpha-\beta}{2}}Q(\|u_i\|_{\mathscr{C}^\alpha},\|u_0\|_{\CC^\alpha})\left(\|u_1-u_2\|_{\mathscr{C}^\alpha}\|u_1^\sharp(t)\|_{\CC^{2\beta-2}}+\|u_1^\sharp(t)-u_2^\sharp(t)\|_{\CC^{2\beta-2}}\right)
$$
and therefore 
$$
\sup_{0<t\leq T}t^{\frac{2\beta-\alpha}{2}}\|\phi_1(t)\|_{\CC^{2\beta-2}}\lesssim T^{\frac{\alpha-\beta}{2}}Q\big(\|u_i\|_{\mathscr{C}^\alpha},\|u_0\|_{\CC^\alpha}\big) \|u_1-u_2\|_{\mathbf{C}^\beta_{\alpha,T}(X)}
$$
for yet another positive constant $Q(\cdots)$ with polynomial growth with respect to its parameters.

-- For $\phi_2$, as the function $\sharp_2$ is locally lipschitz in $(u,u')$ it will not cause any issue in the estimate and only the first two terms need a special treatment as they involve a $u^\sharp$ factor. We only consider the $B$-term as they are both dealt with in the same way. We need to estimate
$$
\left\|\big(g'(u_1)\,\PI\big(B(u_1^\sharp) , \xi\big) - g'(u_2)\,\PI\big(B(u_2^\sharp) , \xi\big)\big)(t)\right\|_{\CC^{\alpha+\beta-2}}
$$
for a positive time $t$. An elementary splitting gives for this quantity the upper bound
\begin{equation*} \begin{split}
\left\| \big(g'(u_1)-g'(u_2)\big)(t) \, \PI\big(B(u_1^\sharp) , \xi\big)(t) \right\|_{\CC^{\alpha+\beta-2}} &+ \big\| \big(g'(u_1)-g'(u_2)\big)(t)\big\|_{\CC^{\beta}} \|u_1^\sharp(t)\|_{\CC^{2\beta}}   \\
&+ \|g'(u_2)(t)\|_{\CC^\beta} \big\|(u_1^\sharp-u_2^\sharp)(t)\big\|_{\CC^{2\beta}}.
\end{split} \end{equation*}
Note that as $g$ is $C^3_b$ one has
$$
\|g'(u_2)(t)\|_{\CC^\beta} \leq C\big(\|u_2\|_{\mathbf{C}^\beta_{\alpha,T}(X)}\big)
$$
and
$$
\big\|\big(g'(u_1)-g'(u_2)\big)(t)\big\|_{\CC^{\beta}} \leq C\big(\|u_i\|_{\mathbf{C}^\beta_{\alpha,T}(X)}\big) \|u_1-u_2\|_{\mathbf{C}^\beta_{\alpha,T}(X)}
$$
for yet again polynomial constants. One therefore has
$$
\sup_{0<t\leq T}t^{\frac{2\beta-\alpha}{2}}\|\phi_2(t)\|_{\CC^{\alpha+\beta-2}} \leq C\big(\|u_i\|_{\mathbf{C}^\beta_{\alpha,T}(X)}\big) \|u_1-u_2\|_{\mathbf{C}^\beta_{\alpha,T}(X)}.
$$
Finally one can apply the modified Schauder estimate from Lemma \ref{SchauderWeight} together with the classical Schauder estimates to recover an estimate on $v_1-v_2$

\begin{equation} \label{EqContraction} \begin{split}
    \|v_1-v_2\|_{\mathbf{C}^\beta_{\alpha,T}(X)}&=\|v_1'-v_2'\|_{\mathscr{C}^\beta} + \big\| v_1^\sharp - v_2^\sharp \big\|_{\mathscr{C}^\alpha} + \sup_{0<t\leq T}t^{\frac{2\beta-\alpha}{2}} \big\|u_1^\sharp(t)-u_2^\sharp(t)\big\|_{\CC^{2\beta}}   \\
    &\leq C\Big(\|u_1\|_{\mathbf{C}^\beta_{\alpha,T}(X)}, \|u_2\|_{\mathbf{C}^\beta_{\alpha,T}(X)}\Big) \, T^{\frac{\alpha-\beta}{2}} \|u_1-u_2\|_{\mathbf{C}^\beta_{\alpha,T}(X)}
\end{split} \end{equation}
with the constant $C(\cdot)$ depending continuously on the enhanced noise and polynomially on $u_i\in \mathbf{C}^\beta_{\alpha,T}(X)$. Here it is crucial that in Lemma \ref{SchauderWeight} the constant depends only on the $\CC^\alpha$ norm of $w$ as we want an estimate on $w=u_0^T$ that does not explode as $T>0$ goes to $0$. To conclude note that 
$$
u := \PT_{\frac{B(g(u_0))}{A(f(u_0))}}X + u_0^T
$$ 
defines an element of $\mathbf{C}^\beta_{\alpha,T}(X)$ so if $\widetilde{u}\in \mathbf{\mathcal{B}}_T(r)$ and 
$$
\max\Big(\|u\|_{\mathbf{C}^\beta_{\alpha,T}(X)},\|\Phi(u)\|_{\mathbf{C}^\beta_{\alpha,T}(X)}\Big) < \frac{r}{2}
$$ 
then 
$$
\big\| \Phi(\widetilde{u}) \big\|_{\mathbf{C}^\beta_{\alpha,T}(X)}\leq \frac{r}{2} + 2rC(r) T^{\frac{\alpha-\beta}{2}}\leq r
$$
if $T$ is chosen small enough, depending on $r$. This means that for $r$ large enough and $T$ small enough, both depending only on the size of the initial value $u_0\in \mcC^\alpha(\bbT^2)$ and the problem data 
$$
\Big(\xi , \PI\big(A(X),\xi\big) , \PI\big(B(X) , \xi\big), f, g\Big),
$$ 
the map $\Phi$ sends the ball $\mathbf{\mathcal{B}}_T(r)$ into itself and is a contraction from \eqref{EqContraction}. As such it has a unique fixed point in that ball and the equation
$$
\partial_tu- A(f(u))\Delta u = B(g(u))\,\xi
$$
with initial condition $u(0) = u_0\in \CC^\alpha(\T^2)$ has a unique paracontrolled solution on a small time interval $[0,T]$.

%----------------------------------------%
\section{Renormalization matters}
\label{SectionRenormalization}
%----------------------------------------%

We construct in Subsection \ref{SUBSEC:EnhancedNoise} the enhanced noise; it was used as a data in the analytic setting of Section \ref{SectionPCFormulation}. The renormalization process used for that purpose comes with a dynamical interpretation of the solution to the fixed point problem for $\Phi$. The renormalized equation is described in Subsection \ref{SUBSEC:RenormalizedEquation}.

\medskip

\subsection{Enhanced noise.}{\label{SUBSEC:EnhancedNoise}}
\hspace{.25cm}Following the expansions and given an arbitrary symbol $a$ of order $s\leq0$ we need to make sense of the resonant product $\PI\big(A(X) , \xi\big)$ as an element of $\CC^{2\alpha-s-2}(\T^2)$. This is done via a renormalization procedure described here. Let $\chi : [0,\infty)\rightarrow [0,1]$ be of class $C^\infty$, equal to $1$ on $[0,1]$ and $0$ on $[2,\infty)$. Set for any distribution $u$ on $\T^2$
$$
R^\varepsilon(u) := \F^{-1}\big(\chi(\varepsilon \vert\cdot\vert)\,\widehat{u}\,\big)
$$
and recall that we regularize the space white noise $\xi$ into a smooth function using that smooth Fourier cut-off function
$$
\xi^\varepsilon := R^\varepsilon(\xi).
$$

\smallskip

\begin{thm}{\label{THM:Renorm}}
For any $a\in\Sigma^s_0$ with associated operator $A$ there exists some deterministic diverging function $c:\varepsilon\in(0,1]\to C^\infty(\T^2)$ such that
$$
\PI\big(A(X^{\varepsilon}) , \xi^{\varepsilon}\big) - c^{\varepsilon}
$$
converges in $L^r\big(\Omega; \CC^{2\alpha-2-s}(\T^2)\big)$ for any $1\leq r<\infty$ to a limit $:\hspace{-0.08cm}\PI(AX,\xi)\hspace{-0.08cm}:$ in $\CC^{2\alpha-2-s}(\T^2)$.
\end{thm}

\smallskip

\begin{Dem}
First note that 
$$
\xi^{\varepsilon} \overset{\CC^{\alpha-2}}{\longrightarrow} \xi\quad\text{and}  \quad  X^{\varepsilon} \overset{\CC^{\alpha}}{\longrightarrow} X
$$
almost surely and in $L^r(\Omega)$ for any $1\leq r<\infty$. This is due to continuity of $(-\Delta)^{-1}$ for $X^{\varepsilon}$, and the hypercontractivity and Besov embedding arguments used below can also be used to prove the convergence of $\xi^{\varepsilon}$. Note that due to the commutator estimate from Proposition 3.7 in \cite{BDFT23} one has $[A,R^\varepsilon](X)\in \CC^{2\alpha-s}(\T^2)$. This brings us back to looking at
$$
\PI\big((A(X))^{\varepsilon} , \xi^{\varepsilon}\big) - c^{\varepsilon}
$$
which is what we look at below. Note that if $\widehat{a}=\widehat{a}^x$ is the Fourier transform of $a$ with respect to $x$, then 
$$
(AX)^{\varepsilon}=\sum_{n\in\Z^2}\chi(\varepsilon|n|)\widehat{AX}(n)e_n=\sum_{n,p\in\Z^2}\chi(\varepsilon|n|)\widehat{X}(n-p)\widehat{a}(p,n-p)e_n =: \sum_{p\in\Z^2}(AX)^{\varepsilon}_p,
$$
where $e_n$ is the function $e_n(x)=e^{in\cdot x}$. We shall carry out the renormalization procedure for each $(AX)^{\varepsilon}_p$ and then pass to the limit by keeping careful track of the $p$-dependency.\\
Fix some $p\in\N$ and $\varepsilon>0$, the very definition of $\PI$ yields
$$
\PI\big((AX)^{\varepsilon}_p,\xi^{\varepsilon}\big) = \sum_{|i-j|\leq1}\sum_{k,l\in\Z^2,l\neq p}\rho_i(k)\,\rho_j(l)\,\chi(\varepsilon|k|)\,\chi(\varepsilon|l|)\,\frac{\widehat{a}(p,l-p)}{|l-p|^2}\,\widehat{\xi}(k)\,\widehat{\xi}(l-p)e_{k+l}.
$$
Set
\begin{align}{\label{EQ:Definitioncp}}
c_p(\varepsilon) := \E\left[\PI\big((AX)^{\varepsilon}_p,\xi^{\varepsilon}\big)\right] = \sum_{|i-j|\leq1}\sum_{k\in\Z^2\setminus\{0\}} \rho_i(k)\,\rho_j(p-k)\,\chi(\varepsilon|k|)\,\chi(\varepsilon|p-k|)\,\frac{\widehat{a}(p,k)}{|k|^2}\,e_p
\end{align}
and 
$$
Y_p^{\varepsilon} := (AX)_p^{\varepsilon} - c_p(\varepsilon).
$$
We prove that $Y_p^{\varepsilon}$ converges towards some $Y_p$ in $\CC^{2\alpha-2}(\T^2)$, in $L^r(\Omega)$ for any $r\geq1$. Let $q\geq-1$ and $\varepsilon,\eta>0$. Then 
\begin{eqnarray*}
\Delta_q(Y^{\varepsilon}_p-Y^\eta_p)&=&\sum_{\substack{|i-j|\leq1\\ k,l\in\Z^2,l\neq p}} \rho_q(k+l)\,\rho_i(k)\,\rho_j(l)\,\frac{\widehat{a}(p,l-p)}{|l-p|^2} \Big(\chi(\varepsilon |k|)\,\chi(\varepsilon|l|) - \chi(\eta |k|)\,\chi(\eta|l|)\Big)\\
&&\qquad\qquad\times\left(\widehat{\xi}(k)\,\widehat{\xi}(l) - \delta_{k,p-l}\right)e_{k+l}.
\end{eqnarray*}
When computing the expectation $\E\big[|\Delta_q(Y^{\varepsilon}_p-Y^\eta_p)|^2\big]$ we will come across a term in the form of 
$$
\E\left[\left(\widehat{\xi}(k)\,\widehat{\xi}(l-p) - \delta_{k,p-l}\right)\left(\overline{\widehat{\xi}(k')}\,\overline{\widehat{\xi}(l'-p)} - \delta_{k',p-l'}\right)\right].
$$
This can be handled using the properties of the kernel of $\widehat{\xi}(\cdot)$ and Isserlis' formula as $\widehat{\xi}(k)$ are gaussian
$$
\E\left[\left(\widehat{\xi}(k)\,\widehat{\xi}(l-p) - \delta_{k,p-l}\right)\left(\overline{\widehat{\xi}(k')}\,\overline{\widehat{\xi}(l'-p) } - \delta_{k',p-l'}\right)\right] = \delta_{k,k'}\,\delta_{l,l'} + \delta_{k,l'-p}\,\delta_{l,k'+p}
$$
hence we can keep only the terms where $(k,l)\in\left\{(k',l'),(l'-p,k'+p)\right\}$. This allows to get rid of the $e_{k+l-k'-l'}$ terms and get
\begin{align*}
&\E\left[\big|\Delta_q(Y^{\varepsilon}_p-Y^\eta_p)\big|^2\right] = \sum_{\substack{|i_1-j_1|\leq1\\ |i_2-j_2|\leq1}}\sum_{\substack{k,l\in\Z^2, l\neq p\\ k',l'\in\Z^2, l'\neq p}}\rho_q(k+l)^2\,\rho_{i_1}(k)\,\rho_{i_2}(k')\,\rho_{j_1}(l)\,\rho_{j_2}(l')   \\
&\hspace{2cm}\times \Big\{\chi(\varepsilon |k|)\,\chi(\varepsilon|l|) \hspace{-0.05cm}-\hspace{-0.05cm} \chi(\eta |k|)\,\chi(\eta|l|)\Big\} \hspace{-0.03cm} \Big\{\chi(\varepsilon |k'|)\chi(\varepsilon|l'|) \hspace{-0.05cm}-\hspace{-0.05cm} \chi(\eta |k'|)\chi(\eta|l'|)\Big\}   \\
&\hspace{2cm}\times \frac{1}{|l-p|^2}\,\frac{1}{|l'-p|^2}\,\widehat{a}(p,l-p)\,\overline{\widehat{a}(p,l'-p)}\big(\delta_{k,k'}\,\delta_{l,l'} + \delta_{k,l'-p}\,\delta_{l,k'+p}\big).
\end{align*}
Since in both sums $|i-j|\leq1$ and $\text{supp}(\rho)$ is a fixed annulus, we can actually keep only the terms where $|k|\sim|l|$ (resp. $|k'|\sim|l'|$), that is $\theta|l|\leq|k|\leq \theta^{-1}|l|$ for some $\theta>0$ that does not depend on $i_1,i_2,j_1,j_2$ and bound the sum on $i_1,i_2,j_1,j_2$ by a constant independent on $k,l,k',l'$, leaving us with an upper bound for $\E\left[\big|\Delta_q(Y_p^{\varepsilon}-Y_p^\eta)(x)\big|^2\right]$ of the form

\begin{align*}
&\sum_{\substack{k,l\in\Z^2, l\neq p\\ \theta|l|\leq|k|\leq \theta^{-1}|l|}} \hspace{-0.4cm}\rho_q(k+l)^2\,\frac{\Psi_1^{\varepsilon,\eta}(k,l)}{|l-p|^4}\,\big|\widehat{a}(p,l-p)\big|^2\\
&\hspace*{3cm}+\sum_{\substack{k,l\in\Z^2, k\neq 0, l\neq p\\ \theta|l|\leq|k|\leq \theta^{-1}|l|}}\rho_q(k+l)^2\,\frac{\Psi_2^{\varepsilon,\eta}(k,l,p)}{|k|^2|l-p|^2}\,\big|\widehat{a}(p,l-p)\widehat{a}(p,k)\big|
\end{align*}
where
$$
\Psi_1^{\varepsilon,\eta}(k,l) = \big|\chi(\varepsilon |k|)\,\chi(\varepsilon|l|) - \chi(\eta |k|)\,\chi(\eta|l|)\big|^2
$$
and
$$
\Psi_2^{\varepsilon,\eta}(k,l,p) = \Big|\chi(\varepsilon |k|)\,\chi(\varepsilon|l|) - \chi(\eta |k|)\,\chi(\eta|l|)\Big| \, \Big|\chi(\varepsilon |l+p|)\,\chi(\varepsilon|k-p|) - \chi(\eta |l+p|)\,\chi(\eta|k-p|)\Big|.
$$

\noindent Apply Estimate \eqref{EQ:BoundSymbol} for $k=l-p$, we have 

$$
\left|\frac{\widehat{a}(p,l-p)}{|l-p|^2}\right|\lesssim \langle p\rangle^{-4}\,\langle l\rangle^{s-2}\,\mathbf{1}_{|p|\leq\mu(1+|l-p|)}
$$

\noindent so that, up to the $\langle p\rangle^{-4}$ factor, $p$ only shows up in $\Psi_2\mathbf{1}_{|p|\leq\mu(1+|l-p|)}$. Note that thanks to the size constraint on $|p|$ and $\chi$ being smooth and bounded, we can easily derive a bound for $\Psi_1$ and $\Psi_2$ (uniformly in $p$) as

$$
|\Psi_1^{\varepsilon,\eta}(k,l)|\vee|\Psi_2^{\varepsilon,\eta}(k,l,p)| \lesssim_\gamma |\varepsilon-\eta|^\gamma \big(|k|^\gamma+|l|^\gamma\big)
$$
for any $\gamma\in [0,2]$. Fix $\delta>0$ to be a small parameter and write $n=k+l$. As $|k|\leq \theta^{-1}|l|$, $|n|\leq (1+\theta^{-1})|l|$, this yields
\begin{eqnarray*}
    \E\left[\big|\Delta_q(Y_p^{\varepsilon}-Y_p^\eta)(x)\big|^2\right]&\lesssim&\frac{|\varepsilon-\eta|^\gamma}{\langle p\rangle^8}\sum_{n\in\Z^2}\rho_q(n)^2\sum_{\substack{k+l=n\\ \theta|l|\leq|k|\leq \theta^{-1}|l|}}\frac{|k|^{\gamma}+|l|^\gamma}{\langle l\rangle^{4-2s}}\\
    &\lesssim&\frac{|\varepsilon-\eta|^\gamma}{\langle p\rangle^8}\sum_{n\in\Z^2}\frac{\rho_q(n)^2}{\langle n\rangle^{2-\delta-2s}}\sum_{\substack{k+l=n\\ \theta|l|\leq|k|\leq \theta^{-1}|l|}}\frac{1}{\langle l\rangle^{2+\delta-\gamma}}\\
    &\lesssim&\frac{|\varepsilon-\eta|^\gamma}{\langle p\rangle^8}\sum_{n\in\Z^2}\frac{\rho_q(n)^2}{\langle n\rangle^{2-\delta-2s}}
\end{eqnarray*}

\noindent where the last inequality is valid as soon as $\delta>\gamma$, say $\gamma=\delta/2$. Then it is only a matter of counting how many $\Z^2$ points belong to a $2^q$-sized annulus to finally obtain the bound
$$
\E\left[\big|\Delta_q(Y_p^{\varepsilon}-Y_p^\eta)(x)\big|^2\right]\lesssim 2^{q(\delta+2s)}|\varepsilon-\eta|^{\delta/2}\,\langle p\rangle^{-8}
$$
for any $\delta>0$. Using Gaussian hypercontractivity, this yields
$$
\E\left[\big|\Delta_q(Y_p^{\varepsilon}-Y_p^\eta)(x)\big|^r\right]\lesssim 2^{q\frac{r}{2}(\delta+2s)}|\varepsilon-\eta|^{r\delta/4}\,\langle p\rangle^{-4r}
$$
so that multiplying both sides by $2^{qr(2\alpha-2-s+2/r)}$ and summing over $q\geq-1$, we recover the $B^{2\alpha-2+2/r}_{r,r}(\T^2)$ norm
$$
\E\Big[\big\|Y_p^{\varepsilon}-Y_p^\eta\big\|_{B^{2\alpha-2-s+2/r}_{r,r}}^r\Big] \lesssim \langle p\rangle^{-4r}\left(\sum_{q\geq-1}2^{qr(2\alpha-2+\delta/2+2/r)}\right)|\varepsilon-\eta|^{r\delta/4}
$$
where the series converges if $r>\frac{2}{1-\alpha}$ and $\delta<2-2\alpha$ for instance. Combined with the continuous embedding of $B^{2\alpha-2-s+2/r}_{r,r}(\T^2)$ into $\CC^{2\alpha-2-s}(\T^2)$ this yields the estimate
\begin{align}{\label{EQ:KolmogrovCriterion}}
\E\left[\left\|Y_p^{\varepsilon}-Y_p^\eta\right\|_{\CC^{2\alpha-2-s}}^r\right]\lesssim \langle p\rangle^{-4r}\,|\varepsilon-\eta|^{r\frac{\delta}{4}}
\end{align}
and proves that $(Y_p^{\varepsilon})$ is a Cauchy, therefore convergent, sequence in $L^r(\Omega; \CC^{2\alpha-2-s}(\T^2))$ for any $r$ large enough -- hence for any $r\geq1$ from hypercontractivity. This proves the existence of the renormalized product $:\hspace{-0.08cm}\PI((AX)_p,\xi)\hspace{-0.08cm}:$ for each $p\in\Z^2$. Note that we kept explicit track of $p$. The dominated convergence theorem ensures that the renormalization tranfers to the whole series $\sum_p Y^{\varepsilon}_p$. This gives the $L^r(\Omega)$ convergence in $\CC^{2\alpha-2-s}(\T^2)$ of the renormalized product for all $1\leq r<\infty$.   %\vspace{0.1cm}
\end{Dem}

\noindent \textbf{\textsf{Remark --}} {\it Note that applying the Kolmogorov continuity criterion on the estimate \eqref{EQ:KolmogrovCriterion} for $r\delta>2$ yields the existence of a modification of the process $\varepsilon\mapsto Y_p^{\varepsilon}$ that is a continuous function. The convergence of the renormalized product then also happens almost surely if we allow for a change of the probability space.   }

\ssk

%%--------------------------------------%%
\subsection{Renormalized equation.}
{\label{SUBSEC:RenormalizedEquation}}

\hspace{.25cm}As for the convergence part of the claim, fix $u_0\in\CC^{\alpha}$, note that the fixed point argument we run holds locally uniformly in the enhanced noise data 
$$
\Xi := \Big(\xi,\PI\big(A(X),\xi\big) , \PI\big(B(X) , \xi\big)\Big),
$$ 
hence so does the fixed point solution we obtain. We shall denote by 
$$
\Xi\mapsto \scrS(\Xi)
$$
the solution mapping corresponding to the paracontrolled equation starting from $u_0$; it is a continuous function of the enhanced noise data $\Xi$. If $\zeta$ is a smooth noise, $\zeta$ has a natural enhancement $Z$. Write $\widetilde{\scrS} : \zeta\mapsto v$ for the solution map corresponding to the locally in time well-posed quasilinear equation
$$
\partial_tv - A(f(v))\Delta u = B(g(v))\zeta
$$
with initial condition $u_0$. The map $\scrS$ extends $\widetilde{\scrS}$ in so far as
$$
\scrS(Z) = \widetilde{\scrS}(\zeta).
$$
Denote by 
$$
\zeta\in C^\infty([0,T]\times\T^2)\mapsto\scrS^{\varepsilon}(\zeta)
$$ 
the solution map corresponding to the locally in time well-posed equation
$$          
\partial_t u^\varepsilon - A\big(f(u)\big) \Delta u = B\big(g(u)\big) \, \zeta + c_a(\varepsilon) \bigg(\frac{B\big(g(u)\big)}{A\big(f(u)\big)}\bigg)^2 f'(u) - c_b(\varepsilon) \frac{B\big(g(u)\big)}{A\big(f(u)\big)} \, g'(u)
$$
with initial condition $u_0$. Then we have the identity
$$
\scrS^{\varepsilon}(\xi^{\varepsilon}) = \scrS\Big(\xi^{\varepsilon},\PI\big(A(X^{\varepsilon}) , \xi^{\varepsilon}\big) - c_a(\varepsilon),\PI\big(B(X^{\varepsilon}) , \xi^{\varepsilon}\big) - c_b(\varepsilon)\Big).
$$
Indeed, in view of Equations \eqref{EQ:ExpansionB} and \eqref{EQ:ExpansionA}, replacing the formal resonant products $\PI(B(X),\xi)$ and $\PI(A(X),\xi)$ therein by their approximations $\PI\big(B(X^{\varepsilon}) , \xi^{\varepsilon}\big) - c_b(\varepsilon)$ and $\PI\big(A(X^{\varepsilon}) , \xi^{\varepsilon}\big) - c_a(\varepsilon)$, respectively, and then regrouping the (well-defined) products together we obtain that Equation \eqref{EqReformulation} driven by the noise data
$$
\Big(\xi^{\varepsilon},\PI\big(A(X^{\varepsilon}) , \xi^{\varepsilon}\big) - c_a(\varepsilon),\PI\big(B(X^{\varepsilon}) , \xi^{\varepsilon}\big) - c_b(\varepsilon)\Big)
$$
writes 
$$          
\partial_t u - A\big(f(u)\big) \Delta u = B\big(g(u)\big) \, \xi^\varepsilon + c_a(\varepsilon) \bigg(\frac{B\big(g(u)\big)}{A\big(f(u)\big)}\bigg)^2 f'(u) - c_b(\varepsilon) \frac{B\big(g(u)\big)}{A\big(f(u)\big)} \, g'(u).
$$
Note that since $\scrS$ depends locally uniformly on the noise data, there exists a random time $T$ that does not depend on $\varepsilon$ for which the equation above starting from $u_0$ is well-posed in $\mathscr{C}^\alpha_T$. Together with the continuity of $\scrS$ and Theorem \ref{THM:Renorm}, this completes the proof of the convergence.

\medskip

\noindent \textbf{\textsf{Remark --}} {\it If we ask further that the action of $A$ and $B$ is diagonal -- if for instance $A$ and $B$ are convolution operators, then a strong consequence of Theorem \ref{THM:Renorm} is that the renormalization function is a constant. It would be interesting to see whether the renormalization function can be constant in some other cases.   }

\appendix
% \appendixpage
\section{Estimates from paracontrolled calculus}{\label{APP:PC}}

We gather in this section the estimates used in the exposition above. Let $\chi$ and $\rho$ be smooth compactly supported radial functions on $\R^d$. Assume $\chi$ is supported in a ball and $\rho$ is supported in an annulus. Set $\rho_j(\cdot) := \rho(2^{-j}\cdot)$ for $j\geq 0$ and $\rho_{-1} := \chi$. One can find $\chi$ and $\rho$ such that 
$$
\sum_{j\geq -1}\rho_j\equiv1\quad\text{and}\quad \rho_j\rho_i\equiv0\ \text{for}\ |i-j|\geq0
$$
The $i$-th Littlewood-Paley projector is defined as $\Delta_jf := \F^{-1}(\rho_j\widehat{f}\,)$. In those terms one has the formal decomposition
\begin{eqnarray*}
uv&=&\sum_{i<j-1} (\Delta_iu)(\Delta_jv) + \sum_{|i-j|\leq1} (\Delta_iu)(\Delta_jv) + \sum_{i>j+1}(\Delta_iu)(\Delta_jv)   \\
    &=:&\PA_uv+\PI(u,v)+\PA_vu
\end{eqnarray*}
The operator $\PA$ is called paraproduct operator and the operator $\PI$ the resonant operator. One can think of $\PA_uv$ as a function or distribution that globally \textit{looks like} $v$, but with lower frequencies modulated by $u$. This decomposition allows for a precise track of regularity with the following continuity result.

\ssk

\begin{prop} \label{PropContinuityParaproduct}
One has
\begin{equation*}
\|\PA_uv\|_{B^{\alpha\wedge0+\beta}_{p,q}}\lesssim\|u\|_{B^{\alpha\wedge 0}_{p_1,q}}\|v\|_{B^\beta_{p_2,q}}
\end{equation*}
for any $\alpha,\beta\in\R$ and $1\le p_1,p_2,p,q\le\infty$ such that $\frac1p=\frac1{p_1}+\frac1{p_2}$, and
\begin{equation*}
\|\PI(u,v)\|_{B^{\alpha+\beta}_{p,q}}\lesssim\|u\|_{B^\alpha_{p_1,q}}\|v\|_{B^\beta_{p_2,q}}
\end{equation*}
for any $\alpha,\beta\in\R$ such that $\alpha+\beta>0$. 
\end{prop}

\ssk

In view of the regularity condition $\alpha+\beta>0$, it is clear that the ill-defined part of a product comes from its resonant term. This can however be circumvented for functions having themselves a paracontrolled expansion. This is the purpose of the following operator, first introduced in \cite{GIP15}, see Lemma 2.4 therein. Define the operator 
$$
\DC : (u,v,w)\mapsto \PI(\PA_uv,w) - u\,\PI(v,w)
$$
on the set of smooth functions on $\T^2$.

\ssk

\begin{prop} \label{Corrector}
Pick some regularity exponents $\alpha\in(0,1)$, $\beta,\gamma\in\R$ such that $\beta+\gamma<0$ and $\alpha+\beta+\gamma>0$. The operator $\DC$ has a unique extension into a bounded trilinear operator from $\CC^\alpha(\T^2)\times \CC^\beta(\T^2)\times \CC^\gamma(\T^2)$ to $\CC^{\alpha+\beta+\gamma}(\T^2)$.
\end{prop}

\ssk

With these results at hand, we could already almost tackle linear singular SPDEs such as the Parabolic Anderson Model (PAM) 
$$
\partial_t u - \Delta u = u\xi
$$
on the 2-dimensional torus $\T^2$. The idea is then to work in a random subspace of functions having a paracontrolled expansion
\begin{equation}{\label{PCExpansion}}
    u=\PA_{u'}X+u^\sharp
\end{equation}
where $X$ depends only on the noise $\xi$ and the remainder $u^\sharp$ has better regularity than $u$ itself. We then run a fixed point argument on the pair $(u',u^\sharp)$ rather than $u$. Note that the corrector $\DC$ is defined to get around singular products, but other commutators and correctors might be needed to investigate well-defined products that may not be written in a good form for the analysis. A thorough study of these commutators in a very general framework can be found in the works of Bailleul \& Bernicot \cite{BB16} and \cite{BB19}. We only retain the following result on the merging operator that can be found in Theorem 3.6 of \cite{BB16}.

\ssk

\begin{prop} \label{PropIteratedParaproducts}
For $0<\beta\leq\alpha<1$ and $h_1\in \CC^\alpha(\T^2), h_2\in C^\beta(\T^2)$ and $h_3\in \CC^\alpha(\T^2)$ one has
$$
\big\Vert \PA_{h_1}\big(\PA_{h_2}h_3\big) - \PA_{h_1h_2}h_3 \big\Vert_{\CC^{\alpha+\beta}} \lesssim \Vert h_1\Vert_{\CC^{\alpha}} \Vert h_2\Vert_{\CC^{\beta}} \Vert h_3\Vert_{\CC^{\alpha}}.
$$
\end{prop}

\ssk

As for the non-linear case, we need to further investigate the action of non-linear functions on paraproducts, this is done via Bony's paralinearization lemma. For $f\in C^2_b(\R)$ and $u\in \CC^\alpha(\T^2)$ with $\alpha>0$ set
$$
R_f(u) := f\circ u - \PA_{f'\circ u}u.
$$
then we have the paralinearization estimate -- see for instance Theorem 5.2.5 in \cite{M08} or Theorem 2.92 in \cite{BCD}.

\ssk

\begin{prop} \label{Paralinearization}
For $f\in C^2_b(\R)$ and $u\in \CC^\alpha(\T^2)$ with $\alpha>0$ one has
    $$
    \left\|R_f(u)\right\|_{\CC^{2\alpha}} \lesssim \|f\|_{C^2_b}(1+\|u\|_{\CC^\alpha}^2)
    $$
Moreover we have a local-Lipschitz estimate
    $$
    \big\|f\circ u-f\circ v\big\|_{\CC^\alpha}\leq \|f\|_{C^2_b}\big(1+\|u\|_{\CC^\alpha}\big) \|u-v\|_{\CC^\alpha}
    $$
\end{prop}

\ssk

This can be further improved as in Lemma C.1. of \cite{GIP15} to the special case of functions having an expansion with a smoother remainder term as in \eqref{PCExpansion}.

\ssk

\begin{prop} \label{PropRefinedParalinearization}
Pick $0<\beta\leq\alpha<1$ such that $\alpha+\beta>1$. For $u\in \CC^\alpha$ and $v\in \CC^{\alpha+\beta}$, for $f\in C^3_b(\R)$, one has
$$
\Big\|f(u+v)-\PA_{f'(u+v)}(u+v)\Big\|_{\CC^{\alpha+\beta}}\lesssim \|f\|_{C^3_b}\big(1+\|u\|^{1+\beta/\alpha}_{\CC^\alpha}+\|v\|^2_{L^\infty}\big)(1+\|v\|_{\CC^{\alpha+\beta}})
$$
\end{prop} 

\ssk

To take into account the time dependency we now define Gubinelli, Imkeller \& Perkowski's modified paraproduct $\PT$ on parabolic functions as in \cite{GIP15}. Fix some smooth function $\varphi:\R\to\R_+$ with compact support and integral $1$. For $f\in C_TC^\infty$ and $i\geq -1$ define for $0\leq t\leq T$
$$
Q_i(f)(t) := \int_\R 2^{2i}\varphi\left(2^{2i}(t-s)\right)f\big((s\wedge T)\vee0 \big)ds
$$
and
$$
\PT_fg := \sum_{i<j-1}\Delta_i\big(Q_j(f)\big) (\Delta_j g).
$$
The operator $\PT$ satisfies the following uniform in time continuity estimates
$$
\left\|\left(\PT_fg\right)(t)\right\|_{\CC^\alpha}\lesssim \|f\|_{C_TL^\infty}\|g(t)\|_{\CC^\alpha}.
$$
Denote by $\L$ the heat operator $\partial_t-\Delta$ and write $\Vert\cdot\Vert_{C_T\CC^\beta}$ for the natural norm on $C\big([0,T],\CC^\beta(\T^2)\big)$. Similarly we write $\Vert\cdot\Vert_{C^{\alpha/2}_TL^\infty}$ for the natural norm on $C^{\alpha/2}\big([0,T],L^\infty(\T^2)\big)$, then the modified paraproduct $\PT$ satisfies the following commutation estimates -- see Lemma 5.1 in \cite{GIP15}.

\ssk

\begin{prop} \label{CommutationParaproducts}
Pick two regularity exponents $\alpha\in(0,1), \beta\in \R$ and a finite positive time horizon $T$. For $u\in\mathscr{C}^\alpha_T$ and $v\in C\big([0,T],\CC^\beta(\T^2)\big)$ one has
\begin{enumerate}
	\item[(a)] $\;\left\|\L\left(\PT_uv\right)-\PT_u(\L(v))\right\|_{C_T \CC^{\alpha+\beta-2}}\lesssim \|u\|_{\mathscr{C}^\alpha_T}\|v\|_{C_T \CC^\beta}$,
        \item[(b)] $\;\left\|\PT_uv - \PA_uv\right\|_{C_T \CC^{\alpha+\beta}}\lesssim \|u\|_{C_T^{\alpha/2} L^\infty}\|v\|_{C_T \CC^\beta}.$
\end{enumerate}
\end{prop}

\ssk

The Schauder estimates for the inverse heat operator $\L^{-1}$ take the following form -- see e.g. Lemma 2.5 in \cite{PV23} and references therin.

\ssk

\begin{prop} \label{SchauderClassic}
The following statements holds true.
    \begin{enumerate}
        \item[(a)] For $t>0$, $\alpha\in\R$ and $\delta\geq 0$ we have
        $$
        \left\|e^{t\Delta}(u)\right\|_{\CC^{\alpha+\delta}}\lesssim t^{-\delta/2}\|u\|_{\CC^\alpha}\quad\text{and}\quad\left\|e^{t\Delta}(u)\right\|_{\CC^{\delta}}\lesssim t^{-\delta/2}\|u\|_{L^\infty}.
        $$
        
        \item[(b)] For $\alpha\in(0,1)$, $u\in \CC^\alpha$, then for any $t\geq0$ one has
        \begin{equation} \label{EqHeatEstimate}
        \left\|e^{t\Delta}(u) - u\right\|_{L^\infty}\lesssim t^{\alpha/2}\|u\|_{\CC^\alpha}.
        \end{equation}
        
        \item[(c)] Let $u\in C\big([0,T],\CC^\alpha(\T^2)\big)$ for some $\alpha\in\R$ and $T>0$, then for any $\gamma\in [0,1)$ and $t\in(0,T]$ we have 
        $$
        t^\gamma\left\|\left(\L^{-1}(u)\right)(t)\right\|_{\CC^{\alpha+2}}\lesssim\sup_{s\in[0,t]}s^\gamma\|u(s)\|_{\CC^\alpha}.
        $$
        Moreover if $\alpha\in(-2,0)$ then one has
        $$
        \left\|\L^{-1}(u)\right\|_{C_T^{(\alpha+2)/2}L^\infty}\lesssim\|u\|_{C_T \CC^\alpha}.
        $$
        
        \item[(d)] For $\alpha\in(0,2)$ we have the $T$-uniform estimate
        $$
        \left\|\L^{-1}(u)\right\|_{\mathscr{C}^\alpha_T}\lesssim (1+T)\|u\|_{C_T\CC^{\alpha-2}}.
        $$

	\item[(e)] For $\alpha\in(0,2)$, uniformly in $\varepsilon\in(0,1)$ and $T$, we have 
    $$
    \left\|\L^{-1}(u)\right\|_{\mathscr{C}^\alpha_T}\lesssim T^\varepsilon\|u\|_{C_T\CC^{\alpha-2+2\varepsilon}}.
    $$

	\item[(f)] For $\delta\in[0,1)$, $\alpha\in\R$ and $\beta\in[-\alpha,2-\alpha)$ one has
    $$
    \sup_{0\leq s\leq t}s^{\delta}\big\|\L^{-1}(w)(s)\big\|_{\CC^{\alpha+\beta}}\lesssim T^{\frac{\alpha-\beta}{2}}\sup_{0\leq s\leq t}s^\delta\|w(s)\|_{\CC^{2\alpha-2}}.
    $$
\end{enumerate}
\end{prop}

\ssk

The following elementary statement is useful to get a small factor and produce contractive maps.

\ssk

\begin{prop} 
    For $0<\beta<\alpha<2$ and $u\in\mathscr{C}^\alpha_T$ that is null at time $0$ one has
    $$
    \|u\|_{\mathscr{C}^\beta_T}\lesssim T^{\frac{\alpha-\beta}{2}}\|u\|_{\mathscr{C}^\alpha_T}.
    $$
\end{prop}

\smallskip

Finally the following commutation lemma is useful whenever we deal with a weighted heat operator. See Lemma 2 in \cite{BDH16} for a proof.

\smallskip

\begin{lem}
\label{CommutationWeightedLaplacian}
    For $v\in\mathscr{C}^\beta_T$ and $w\in \CC^{2\beta}(\T^2)$ one has
    $$
    \left\|(\partial_t-w\Delta)\left(\PT_vX\right)-\PA_{wv}\big(-\Delta X)\right\|_{C_T\CC^{\alpha+\beta-2}}\lesssim \big(1+\|w\|_{\CC^{2\beta}}\big) \|v\|_{\mathscr{C}^\beta_T} \|X\|_{\CC^\alpha}.
    $$
\end{lem}

\vspace{0.7cm}

\noindent \textcolor{gray}{$\bullet$} I. Bailleul -- Univ Brest, CNRS, LMBA - UMR 6205, F- 29238 Brest, France.\\
{\it E-mail}: ismael.bailleul@univ-brest.fr

\medskip

\noindent \textcolor{gray}{$\bullet$} H. Eulry -- Univ Rennes, CNRS, IRMAR - UMR 6625, F- 35000 Rennes, France.\\
{\it E-mail}: hugo.eulry@ens-rennes.fr

\end{document}